\documentclass[12pt]{amsart}

\usepackage{amsmath,amssymb,enumerate}

\usepackage{epsfig,fancyhdr,color}

\usepackage{amssymb}
\usepackage{amsmath,amsthm}
\usepackage{bbm}
\usepackage{latexsym}
\usepackage{amscd}
\usepackage{psfrag}
\usepackage{graphicx}
\usepackage[all]{xy}

\usepackage{rotating}
\usepackage{adjustbox}

\usepackage[font=small,labelfont=bf]{caption}

\usepackage{mathtools}

\usepackage{tikz-cd}

\usepackage{stackrel}

\usepackage[utf8]{inputenc}
\usepackage{amsmath}
\usepackage{amssymb}

\usepackage{mathabx,epsfig}
\def\acts{\mathrel{\reflectbox{$\righttoleftarrow$}}}

\usepackage[utf8]{inputenc}

\numberwithin{equation}{section}

\definecolor{NoteColor}{rgb}{1,0,0}

\newtheorem{theorem}{\rm\bf Theorem}[section]

\newtheorem{lemma}[theorem]{\rm\bf Lemma}
\newtheorem{corollary}[theorem]{\rm\bf Corollary}
\newtheorem*{theorem 1}{\rm\bf Proposition 1}
\newtheorem*{theorem 2}{\rm\bf Proposition 2}

\theoremstyle{definition}

\theoremstyle{remark}

\def\interieur#1{\mathord{\mathop{\kern 0pt #1}\limits^\circ}}

\def\interieur#1{\mathord{\mathop{\kern 0pt #1}\limits^\circ}}

\def\hyperp{{\rm I}\kern-.3ex{\rm H}}

\def\hyperl{\mathbb H}

\usepackage{mathabx,epsfig}
\def\acts{\mathrel{\reflectbox{$\righttoleftarrow$}}}

\begin{document}

\title[Angle Defect for Super Triangles]{Angle Defect for Super Triangles}
\
\author{Robert Penner}
\address {\hskip -2.5ex Institut des Hautes \'Etudes Scientifiques\\
35 route des Chartres\\
Le Bois Marie\\
91440 Bures-sur-Yvette\\
France\\
{\rm and}~Mathematics Department,
UCLA\\
Los Angeles, CA 90095\\USA}
\email{rpenner{\char'100}ihes.fr}

\thanks{It is a pleasure to thank Anton Zeitlin for input and guidance several years ago, when we abandoned our collaborative attempt
prove the theorem, and both Athanase Papdopoulos and Dennis Sullivan for more recent enlightening remarks.}

 \date{\today}


\begin{abstract}
We prove that the angle defect minus the area of a super hyperbolic triangle is not identically zero and explicitly compute
the purely fermionic difference.  This disproves the Angle Defect Theorem for ${\mathcal N}=1$ super hyperbolic geometry and
provides a novel non-trivial additive function of super triangles.  The proof techniques involve the orthosymplectic
group ${\rm OSp}(1|2)$ in its action on the super Minkowski space ${\mathbb R}^{2,1|2}$ and brute-force computation.
\end{abstract}

\maketitle





{\let\thefootnote\relax}
\setcounter{footnote}{0} 

\section*{Introduction}
This paper continues our ongoing study of triangles in ${\mathcal N}=1$ super hyperbolic geometry, namely, the geometry of the usual 
two-real dimensional hyperbolic plane enriched with two real anti-commuting variables.  
We proved the Law of Cosines for super hyperbolic triangles in  \cite{norbert} as a birthday present for Norbert A'Campo, only
to find exactly the classical formula with no corrections; however, we also found that this classical formula
relates expressions which do involve so-called {\sl fermionic corrections}.  Here as a birthday present for Athanase Papadopoulos,
we compute the non-zero correction to the Angle Defect Theorem for super hyperbolic triangles.

The classical hyperbolic Angle Defect Theorem states most
elegantly that the sum $S$ of the interior angles of a hyperbolic triangle plus its area $A$ is constant equal to $\pi$.  That is,
$S\geq 0$ is always at most $\pi$, the {\sl defect} $D$ is $D=\pi-S\geq 0$, and the theorem is that $A=D$.  

It was Lambert \cite{lambert} who originally proved this result in hyperbolic geometry in 1766 in a manuscript unpublished in his
lifetime, which finally appeared in German \cite{German} in 1895, after Johann Bernoulli III, who was in charge
of Lambert's estate, considered publishing it
and decided against it.  A French translation \cite{French} by Athanase and his student Guillaume Th\`eret appeared only recently in 2014. 
 
As was the fashion for centuries, and indeed for millennia, to prove the impossibility of the dubious hyperbolic geometry,
theorems were proved within its confines in the hopes of deriving a logical contradiction.  This was Lambert's
motivation, and his elegant proof comprises all of two statements:  Angle defect is additive for triangles, and 
area is the only such additive function, at least up to an overall scalar.  The first assertion is obvious, as is the additivity of 
area, but the uniqueness part of Lambert's second assertion seems to anticipate developments in measure theory by several centuries.

Contemporary textbook proofs of the Angle Defect Theorem typically involve first computing the area of a
triangle with one or more vanishing interior angles and again rely on additivity of area.
This approach is not available in super hyperbolic geometry due to the divergence of the fermionic correction to the area
of a super hyperbolic triangle with vanishing interior angles, as we prove in the Appendix.

Lambert's proof is likewise unavailable in our context because {\sl area is not the only additive function
for super hyperbolic triangles}.  This follows from our main result Theorem \ref{mainthm}, which explicitly computes the non-zero fermionic
correction to the Angle Defect Theorem for super triangles.  This furthermore raises the intriguing question of what are
the additive functions on super hyperbolic triangles.

The classical Angle Defect Theorem
is evidently a direct consequence of the Gauss-Bonnet Theorem in constant curvature.
Differential geometric approaches to super geometry are given in \cite{bryce,kessler,Rogers}, and perhaps
an argument different from ours could be based on that.

This entire paper is devoted to the computation of the area and angle defect of a super hyperbolic triangle.  We work in the super Minkowski space
${\mathbb R}^{2,1|2}$, which is recalled and discussed in Section \ref{sec:mink}, taking the upper sheet $\hyperp$ of the unit hyperboloid
as our model for the super hyperbolic plane.  As has been noted elsewhere \cite{PZ}, this is the most natural model since in the super analogues
of the upper half plane and Poincar\'e disk, not every geodesic ray is asymptotic to a point of the absolute; in contrast, every super geodesic in $\hyperp$ is indeed
asymptotic in each direction to a ray in the super light cone, as we shall recall.

Super geodesics are discussed in Section \ref{sec:geod}.  The real orthosymplectic group ${\rm OSp}(1|2)$, which plays the role here of
${\rm PSL}(2,{\mathbb R})$ or ${\rm SL}(2,{\mathbb R})$ in the classical case, is recalled and explained in Section \ref{sec:ortho}.
The ${\rm OSp}(1|2)$-invariant area form $\Omega$ on $\hyperp$ is computed in Section \ref{sec:omega} as the pull-back of the usual
${\rm PSL}(2,{\mathbb R})$-invariant area form in the classical case, and a primitive $\omega$ with ${\rm d}\omega=\Omega$ is
derived.  Our approach to computing areas is given by Stokes' theorem and hence devolves to computing line integrals of $\omega$
on the edges of a super triangle; additivity of area is thus obvious with this formalism.  

These calculations are at some points concise and rewarding and at other points routine and tedious.
They comprise the bulk of the paper in Sections \ref{sec:line} and \ref{sec:area}.  The final result provides a purely fermionic additive function on super triangles,
which is given by the angle defect minus the area.  It is obviously not identically zero, but it is complicated and far from illuminating.  
With a large
fermionic correction to the classical formula, this is therefore diametrically
opposite to what we found for the Law of Cosines in the earlier work.

Let us emphasize that whereas the classical results we extend involve two-real dimensional geometry,
their super analogues include also two real fermions and hence are higher dimensional.  For instance,
Theorem 6.2 in \cite{norbert} provides the further conditions to guarantee the intersection in $\hyperp$ of two super geodesics
whose underlying classical geodesics intersect in a point.
Exactly what, if anything, our findings might suggest about 
classical or super scissors congruence in dimensions 2, 3 or 4
remains to be seen.  
It seems feasible to follow the scheme here and find an analogous correction for the volume of a super hyperbolic 3-simplex.

This paper is dedicated to Athanase Papadopoulos on the happy occasion  of his 65th birthday.
He is my dearest friend.  We met as  post-docs, and since our doctoral theses had substantial overlap, we presumably
should have been natural competitors.  It was quite the contrary, and from the moment we first shook hands 
at a conference in Warwick, we have been fast friends.

I had just moved from Princeton to Los Angeles, and he from Orsay to Strasbourg.  I more or less immediately joined Athanase in
Strasbourg for a semester, and he directly thereafter joined me in Los Angeles for a semester, each of us with our
young families in tow.  We wrote two papers \cite{PP1,PP2} together at that time on Thurston theory, specifically on pseudo-Ansov mappings 
and measured foliations.  I appreciated Strasbourg more than he enjoyed Los Angeles, and I have often been back to visit him
over the last 40 years.  Strasbourg has been a second home to me, or perhaps more precisely third.

We had three substantial collaborations after that period:    In \cite{PP3,PP4}, we showed that the Weil-Petersson K\"ahler form
on Teichm\"uller space extends naturally to the Thurston two-form on the space of (projective) measured foliations, an early paradigm for tropical geometry as well as for the 
theory of cluster varieties, which co-opted some of our constructions; in \cite{PP5}, we studied so-called broken hyperbolic structures on surfaces,
which have certain controlled discontinuities along geodesics; and in \cite{PP6}, we extended the Hatcher-Thurston theory, among other things, 
to the setting of non-orientable surfaces, providing moves on pants decompositions of a surface so that  finite compositions act transitively.  We have long planned to compute the relations among our 
moves, a project for which I hope we both still have the time, energy and inclination.

We have been through good times and bad times together.  As we have taken our separate trajectories through mathematics and through life from our 
once-similar initial conditions, it has been rewarding to see Athanase develop into the great scholar that he has become.
It has furthermore been an abiding joy in my career and  in my existence to have counted
him among my life-long friends.

\section{Quadrics in Super Minkowski Space}\label{sec:mink}

Let $\hat{{\mathbb R}}=\hat{\mathbb R}[0]\otimes\hat{\mathbb R}[1]$ be the ${\mathbb Z}$/2-graded module over ${\mathbb R}$ with one central generator
$1\in {\mathbb R}\subseteq \hat{\mathbb R}[0]$ of degree zero and countably infinitely many anti-commuting generators $\theta_1,\theta_2,\ldots\in \hat{\mathbb R}[1]$ of degree one.  An arbitrary $a\in\hat{\mathbb R}$ can be written uniquely as
$$a=a_{\#}+\sum_i a_i\theta_i+\sum_{i<j} a_{ij} \theta_i \theta_j++\sum_{i<j<k} a_{ijk} \theta_i \theta_j \theta_k\cdots,$$
where $a_\#, a_i,a_{ij},a_{ijk}\ldots\in {\mathbb R}$.  The  term $a_\#$ of degree zero is called the {\it body} of the {\it super number} $a\in \hat{\mathbb R}$.
If $a\in\hat{\mathbb R}[0]$, then it is said to be an {\it even} super number or  {\it boson}, while if $a\in\hat{\mathbb R}[1]$, then it is said to  be an {\it odd} super number or {\it fermion}. We adopt the notation throughout that fermions are denoted by lower-case Greek letters.

One allows only finitely many anti-commuting factors in any product and only finitely many summands in any super number, a constraint we shall call {\it regularity}.  Regularity implies nilpotence, i.e., for each $a\in\hat{\mathbb R}$, there is some $n\in{\mathbb Z}_{\geq 0}$ with $(a-a_\#)^{n}=0$.  It follows that if $a_\#\neq 0$, then we may write
$$\begin{aligned}
{1\over a}&= {1\over{a_\#+(a-a_\#)}}
={1\over a_\#}~~ {1\over{1+{{a-a_\#}\over a_\#}}}\\
&={1\over a_\#}\biggl [1-{{a-a_\#}\over a_\#}+\bigg({{a-a_\#}\over a_\#}\biggr)^2-\cdots +(-1)^{n-1}\biggl({{a-a_\#}\over a_\#}\biggr)^{n-1}
\biggr ],\\
\end{aligned}$$
and hence $a\in\hat{\mathbb R}$ is invertible if and only if $a_\#\neq 0$.  It similarly follows from regularity that the zero divisors in 
$\hat{\mathbb R}$ are given by the ideal generated by $\hat{\mathbb R}[1]$.  

One analogously extends real-analytic functions
on ${\mathbb R}$ to $\hat{\mathbb R}$ with Taylor series under appropriate restrictions on the body.  For instance for later tacit application, if $b^2=0$ for some $b\in\hat{\mathbb R}$, then for $a\in\hat{\mathbb R}$ with $a_\#\neq 0$, we have
${1\over{a+b}}={1\over a}(1-{b\over {a}})$,
and for $a_\#>0$, we have 
$\sqrt{a+b}=\sqrt{a}~(1+{b\over {2a}})$.

The usual order relation $\leq$ on ${\mathbb R}$ induces one on ${\hat{\mathbb R}}$ with $a\leq b$, for $a,b\in\hat{\mathbb R}$, if and only if $a_\#\leq b_\#$.  Likewise if $a_\#\neq 0$, then the sign $sign(a)$ is defined to be the sign of $a_\#$ and the absolute value is $|a|=sign(a)a$.

{\it Affine $n|m$ dimensional $\hat{\mathbb R}$ super space} is defined to be
$${\mathbb R}^{n\vert m} = \{(x_1,x_2,\ldots,x_n~\vert~\theta_1,\theta_2,\ldots ,\theta_m)\in \hat{\mathbb R}^{n+m}: x_i\in\hat{{\mathbb R}}[0], \theta_j\in \hat{{\mathbb R}}[1]\}.$$
One defines n$\vert$m {\it super manifolds} with charts based on affine super space 
${\mathbb R}^{n|m}$ in the usual way \cite{bryce, kessler, Rogers}, and a {\it (Riemannian) super metric} on an $n\vert m$ super manifold is defined to be a  positive definite boson-valued quadratic form on each tangent space ${\mathbb R}^{n|m}$ as usual.


The principal example for us here is 
{\it (real) super Minkowski {\rm 2,1$\vert$2} space}
$${\mathbb R}^{2,1|2}=\{ (x_1,x_2,y~\vert~\phi,\psi)\in\hat{\mathbb R}^5:x_1,x_2,y\in\hat{\mathbb R}[0]~{\rm and}~ \phi,\psi\in\hat{\mathbb R}[1]\},$$
which supports the boson-valued symmetric bilinear pairing
$$\langle (x_1,x_2,y~\vert~\phi,\psi),(x_1',x_2',y'~\vert~\phi',\psi')\rangle={1\over 2}(x_1x_2'+x_1'x_2)-yy'+\phi\psi'+\phi'\psi$$
with associated quadratic form $x_1x_2-y^2+2\phi\psi$.

The body of ${\mathbb R}^{2,1|2}$ with this inner product is evidently
the classical Minkowski space ${\mathbb R}^{2,1}$ with its (negative) definite restriction to
$$\begin{aligned}
\hyperp'&=\{ {\bf x}=(x_1,x_2,y)\in{\mathbb R}^{2,1}:\langle{\bf x},{\bf x}\rangle=1~{\rm and}~x_1+x_2>0\}\\
\end{aligned}$$
providing a model of the hyperbolic plane, and
we analogously define the {\it super hyperbolic plane} to be
$$\begin{aligned} 
{\hyperp}&=\{ {\bf x}=(x_1,x_2,y~\vert~\phi,\psi)\in{\mathbb R}^{2,1|2}:\langle{\bf x},{\bf x}\rangle=1~{\rm and}~x_1+x_2>0\}\supseteq\hyperp'\\
\end{aligned}$$ 
with its metric likewise induced from the inner product.  Intermediate between the classical and super hyperbolic planes is the {\it bosonic hyperboloid}
$$\widehat{\hyperp}=\{ {\bf x}=(x_1,x_2,y|0,0)\in{\mathbb R}^{2,1|2}:\langle {\bf x},{\bf x}\rangle =1 ~{\rm and}~ x_1+x_2>0\},$$ whose coordinates lie in $\hat{\mathbb R}[0]$, in contrast to those of $\hyperp'$ lying in ${\mathbb R}$. There are the natural inclusions $\hyperp'\subset\widehat{\hyperp}\subset\hyperp$.

Continuing by analogy,
let $$\hyperl'=\{ {\bf h} \in{\mathbb R}^{2,1}:\langle{\bf h},{\bf h}\rangle=-1\}$$ denote the usual hyperboloid of one sheet and  $${\hyperl}=\{ {\bf h} \in{\mathbb R}^{2,1|2}:\langle{\bf h},{\bf h}\rangle=-1\}\supseteq\hyperl'$$ its super analogue.  
Let $L$ denote the collection of isotropic vectors in ${\mathbb R}^{2,1}$ with
$$L^+=\{ {\bf u}=(u_1,u_2,v)\in{L}:~u_1+u_2>0\}$$ 
 the {\it (open) positive light cone} (whose points are affine duals to horocycles in $\hyperp'$, as in \cite{pennerbook}), and let 
${\mathcal L}$ denote the collection of isotropic vectors in ${\mathbb R}^{2,1|2}$ with
$${\mathcal L}^+=\{ {\bf u}=(u_1,u_2,v~\vert~\phi,\psi)\in{\mathcal L}:u_1+u_2>0\}.$$

\section{Super Geodesics}\label{sec:geod}

Super geodesics in $\hyperp$ are parametrized in analogy to the geodesics of general relativity:

\begin{theorem}[Theorem 1.2 of \cite{HPZ}]\label{uvthm}
The general form of a super geodesic  ${\bf x}(s)$ in $\hyperp$ parametrized by arc length $s$ is given by
$$
{\bf x}=\cosh s\,{\bf u}+\sinh s\,{\bf v},
$$
for some ${\bf u}\in\hyperp$, ${\bf v}\in\hyperl$ with $\langle{\bf u}, {\bf v}\rangle =0$. 
The asymptotes of the corresponding super geodesic are given by the rays in ${\mathcal L}^+$
containing the vectors
$\bf e=\bf u+ \bf v$, $\bf f=\bf u-\bf v$.
Conversely,  points ${\bf e},{\bf f}\in{\mathcal L}^+$ with $\langle {\bf e},{\bf f}\rangle =2$ determine a unique corresponding super geodesic, where ${\bf u}={1\over 2}({\bf e}+{\bf f})\in \hyperp$ and ${\bf v}={1\over 2}({\bf e}-{\bf f})\in\hyperl$.   
\end{theorem}

\medskip

\noindent For any ${\bf u}\in\hyperp$ and ${\bf v}\in\hyperl$ with $\langle {\bf u},{\bf v}\rangle =0$, let 
$$L_{{\bf u},{\bf v}}=\{{\rm cosh}\, s~{\bf u}+{\rm sinh}\, s~{\bf v}:s\in{\mathbb R}\}$$ denote the corresponding super geodesic.

\medskip

\begin{proof}
This follows directly from the variational principle applied to the functional 
$$
\int \Big(\sqrt{ |\langle \dot{\bf x}, \dot{\bf x}\rangle|}+\lambda( \langle {\bf x}, {\bf x}\rangle-1)\Big)dt ,
$$
where the dot stands for the derivative with respect to the parameter $s$ along the curve, with
corresponding Euler-Lagrange equations 
$$
{\ddot{\bf x}}=2\lambda {\bf x}, \quad \langle {\bf x}, {\bf x}\rangle=1 \label{geod}
$$
with $s$ is chosen so that $|\langle\dot{\bf x}, \dot{\bf x}\rangle|=1$.
Differentiating two times the second equation and combining with the first equation we have
$$
\lambda=-\frac{\langle\dot{\bf x}, \dot{\bf x}\rangle}{2}.
$$
The solution with $\lambda=-1/2$ can be ruled out, and in the remaining case $\lambda=1/2$, it is expressed as
$$
{\bf x}=\cosh s\, {\bf u}+\sinh s\, {\bf v}.
$$
The constraints $-\langle\dot{\bf x}, \dot{\bf x}\rangle=\langle{\bf x}, {\bf x}\rangle=1$ imply the conditions on ${\bf u}, {\bf v}$. 
\end{proof}

\begin{corollary} \label{geocor}
For any two distinct points ${\bf P},{\bf Q}\in\hyperp$, there is a unique super geodesic containing them,
and the distance $D$ between them satisfies ${\rm cosh}\,D=\langle {\bf P},{\bf Q}\rangle$.
\end{corollary}

\begin{proof}
For existence, note that $\langle{\bf P},{\bf Q}\rangle~>1$ for distinct
${\bf P},{\bf Q}\in\hyperp$, and define
$$\ell=\sqrt{{\langle{\bf P},{\bf Q}\rangle+1}\over{\langle{\bf P},{\bf Q}\rangle-1}},$$
so that $\langle{\bf P},{\bf Q}\rangle={{\ell^2+1}\over{\ell^2-1}}$.  The identity
${\rm cosh}^{-1}t={\rm log}_e(t+\sqrt{t^2-1})$ therefore gives 
$$e^{{\rm cosh}^{-1}\langle{\bf P},{\bf Q}\rangle}=\langle{\bf P},{\bf Q}\rangle+\sqrt{\langle{\bf P},{\bf Q}\rangle^2-1}={{\ell+1}\over{\ell -1}}.$$

We exhibit $L_{{\bf u},{\bf v}}=L_{{{{\bf e}+{\bf f}}\over 2},{{{\bf e}-{\bf f}}\over 2}}$ containing ${\bf P},{\bf Q}$ parametrized as
$${\bf x}(s)={\rm cosh}\, s~{\bf u}+{\rm sinh}\, s~{\bf v}={1\over 2}[{\rm exp}(s)~{\bf e}+{\rm exp}(-s)~{\bf f}]$$
by taking
$$\begin{aligned}
{\bf e}&={{\ell-1}\over{2\ell}}[(1-\ell){\bf P}+(1+\ell){\bf Q}],\\
{\bf f}&={{\ell+1}\over{2\ell}}[(1+\ell){\bf P}+(1-\ell){\bf Q}].\\
\end{aligned}$$
Direct computation confirms that ${\bf e},{\bf f}\in{\mathcal L}^+$ with $\langle{\bf e},{\bf f}\rangle=2$ and
$${\bf x}(0)={\bf P}~{\rm and}~{\bf x}({\rm cosh}^{-1}\langle{\bf P},{\bf Q}\rangle)={\bf Q},$$
thus establishing existence as well as that ${\rm cosh}\,D=\langle {\bf P},{\bf Q}\rangle$. 

\bigskip

For uniqueness, the system of equations
$$\begin{aligned}
{\bf P}&={\rm cosh}\,s~{\bf u}+{\rm sinh}\, s~{\bf v},\\
{\bf Q}&={\rm cosh}\,t~{\bf u}+{\rm sinh}\, t~{\bf v}
\end{aligned}$$
is tantamount to the linear system
$$\begin{pmatrix}{\rm cosh}\, s\,I&{\rm sinh}\, s\,I\\{\rm cosh}\, t\,I&{\rm sinh}\, t\,I\end{pmatrix}~\begin{pmatrix}{\bf u}\\{\bf v}\end{pmatrix}=\begin{pmatrix}{\bf P}\\{\bf Q}\end{pmatrix},$$
where $I$ is the 5-by-5 identity matrix, which has a unique solution since the determinant is non-zero for distinct $s,t$.\end{proof}

\begin{corollary}\label{cor:tgt}
Given distinct ${\bf P},{\bf Q}\in\hyperp$, the unit tangent vector at ${\bf P}$ to the line from ${\bf P}$ to ${\bf Q}$ is given by
 ${{\bf Q}-{\bf P}\,\langle {\bf P},{\bf Q}\rangle}\over{\sqrt{\langle {\bf P},{\bf Q}\rangle^2-1}}$.~\hfill \qedsymbol
 \end{corollary}

It follows from Corollary \ref{geocor} that any three non-collinear points of $\hyperp$ define a {\it super triangle}, namely, three super geodesic segments with disjoint interiors meeting pairwise
at the given points.  The usual Hyperbolic Law of Cosines holds for super triangles, as proved in \cite{norbert} and as follows directly from Corollary \ref{cor:tgt}
upon computing the cosine of an interior angle based on the usual formula.


\section{Orthosymplectic Group ${\rm OSp(1|2)}$}\label{sec:ortho}

We next include basic information concerning the orthosymplectic group ${\rm OSp(1|2)}$, which is among the simplest of Lie super groups and whose body is the classical special linear group ${\rm SL}(2,{\mathbb R})$.  We refer the interested reader to \cite{Kac,Manin}
for more information about general Lie super algebras and super groups
and to \cite{PZ} for details about OSp(1$\vert$2). 

An element $g\in{\rm OSp(1|2)}$ can be represented by the 3-by-3 matrix $g=\begin{psmallmatrix}
a&b&\alpha\\c&d&\beta\\
\gamma&\delta&f\\\end{psmallmatrix},$ where $a,b,c,d,f$ are even and $\alpha,\beta,\gamma,\delta$ are odd,  with multiplication defined by
$$\biggl ( \begin{smallmatrix}
a_1&b_1&\alpha_1\\c_1&d_1&\beta_1\\\gamma_1&\delta_1&f_1\\
\end{smallmatrix}\biggr )
\biggl ( \begin{smallmatrix}
a_2&b_2&\alpha_2\\c_2&d_2&\beta_2\\\gamma_2&\delta_2&f_2\\
\end{smallmatrix}\biggr )=
\biggl (\begin{smallmatrix}
a_1a_2+b_1c_2-\alpha_1\gamma_2&a_1b_2+b_1d_2-\alpha_1\delta_2&a_1\alpha_2+b_1\beta_2+\alpha_1 f_2\\
c_1a_2+d_1c_2-\beta_1\gamma_2&c_1b_2+d_1d_2-\beta_1\delta_2&c_1\alpha_2+d_1\beta_2+\beta_1 f_2\\
\gamma_1 a_2+\delta_1 c_2+\delta_1\gamma_2&\gamma_1 b_2+\delta_1 d_2+f_1 \delta_2&-\gamma_1\alpha_2-\delta_1\beta_2+f_1f_2\\
\end{smallmatrix}\biggr ),$$
where $g$ is required to satisfy the two further conditions: 

$\bullet$~the {\it super determinant} or {\it Berezinian} of $g$ is unity, namely,
$$
{\rm sdet}~g~=~f^{-1}
~\det\left[
\begin{pmatrix}
a&b\\ c&d
\end{pmatrix}
+f^{-1}
\begin{pmatrix}
\alpha\gamma&\alpha\delta\\ \beta\gamma&\beta\delta
\end{pmatrix}
\right]=1;
$$ 
and 

$\bullet$~$g$ is {\it orthosymplectic}, namely
$$g^{st}Jg=J,$$
where
$
J= \begin{psmallmatrix}
0 & -1 &0 \\
1 &\hskip1.2ex 0 &0 \\
0 &\hskip1.2ex 0 & 1 \end{psmallmatrix}
$
and the {\it super transpose} $g^{st}=\begin{psmallmatrix}
~a & ~c & \gamma \\
~b & ~d & \delta \\
-\alpha & -\beta & f \\
\end{psmallmatrix}.$

\medskip

The extra minus signs in the definition of matrix multiplication come from fermionic anti-commutation 
and our specification of the ordering of rows and columns.  The Berezinian
is the analogue of the classical determinant \cite{Kac,Manin,Rogers} and is characterized by being a  multiplicative homomorphism satisfying sdet exp $g$ = exp$(a+d-f)$, but unlike the classical determinant, it is only defined for invertible super matrices.  The functional equation $g^{st}Jg=J$ explains the analogy with the classical symplectic group and provides the simple expression
$g^{-1}=J^{-1}g^{st}J=\begin{psmallmatrix}
~~d & -b & ~~~~\delta \\
-c & ~~a & -\gamma \\
-\beta & ~~\alpha & ~~f \end{psmallmatrix}$
for inversion in ${\rm OSp}(1|2)$.  More explicitly, the functional equation is equivalent to the contraints
\begin{eqnarray}\nonumber
&&\alpha=b\gamma-a\delta, \quad \beta =d\gamma-c\delta, \quad \hskip 2.8ex f=1+\alpha\beta,\nonumber\\
&&\gamma=a\beta-c\alpha, \quad \delta=b\beta-d\alpha, \quad f^{-1}=ad-bc.\nonumber
\end{eqnarray}

\medskip

Notice that there is a bosonic Lie super group
$$\widehat{\rm SL}_2={\rm SL}(2,\hat{\mathbb R}[0])=\{\begin{psmallmatrix}a&b&0\\c&d&0\\0&0&1\\
\end{psmallmatrix}:a,b,c,d\in\hat{\mathbb R}[0]~{\rm and}~ad-bc=1\}$$ together with canonical inclusions
${\rm SL}(2,{\mathbb R})<\widehat{\rm SL}_2<{\rm OSp}(1|2)$.
Another canonical subspace $\{ u(\alpha,\beta):\alpha,\beta\in\hat{\mathbb R}[1]\}\subset{\rm OSp}(1|2)$, which is not a subgroup, is parametrized by
a pair of fermions $\alpha,\beta$ and defined by
$$u(\alpha,\beta)=\begin{pmatrix}1-{{\alpha\beta}\over 2}&0&\alpha\\0&1-{{\alpha\beta}\over 2}&\beta\\\beta&-\alpha&1+\alpha\beta\\\end{pmatrix}\in{\rm OSp}(1|2).$$

\medskip
 
 It is not difficult to prove the following useful

\begin{lemma}\label{factor}
Any element $g\in{\rm OSp}(1|2)$ can be written uniquely as a product
$$g=\begin{psmallmatrix}a&b&0\\c&d&0\\0&0&1\\\end{psmallmatrix}\,u(\alpha,\beta)=u(a\alpha+b\beta,c\alpha+d\beta)\, \begin{psmallmatrix}a&b&0\\c&d&0\\0&0&1\\\end{psmallmatrix},$$
in ${\rm OSp}(1|2)$, for some $\begin{psmallmatrix}a&b&0\\c&d&0\\0&0&1\\\end{psmallmatrix}\in \widehat{\rm SL}_2$ and fermions $\alpha,\beta\in\hat{\mathbb R}[1]$.\hfill\qedsymbol
\end{lemma}

\medskip

Minkowski three-space ${\mathbb R}^{2,1}\approx {\mathbb R}^3$,
as the space of binary symmetric bilinear forms, is naturally coordinatized by 
$$A=\begin{pmatrix}z-x&y\\y&z+x\end{pmatrix}=\begin{pmatrix}x_1&y\\y&x_2\end{pmatrix}\in{\mathbb R}^{2|1},$$
and $g$ in the component  $ {\rm SO}_+(1,2)$ of the identity in $ {\rm SO}(1,2)$ acts via orientation-preserving isometry on $A\in{\mathbb R}^{2,1}$ as change of basis via the adjoint
$$g:A\mapsto g.A=g^tAg.$$  This action of ${\rm SL}(2,{\mathbb R})$ or ${\rm SO}_+(1,2)\approx {\rm PSL}(2,{\mathbb R})$ as the group of isometries of $\hyperp'$ links hyperbolic geometry and elementary number theory.

\bigskip

Likewise, ${\mathbb R}^{2,1|2}\approx {\mathbb R}^{3|2}$ is naturally coordinatized
by $$A=\begin{pmatrix}~x_1&~~y&\phi\\~~y&~x_2&\psi\\-\phi&-\psi&0\\\end{pmatrix}\in{\mathbb R}^{2,1|2},$$
and $g\in{\rm OSp}(1\vert 2)$ acts on $A$  as change of basis again via the adjoint
$$g:A\mapsto g.A=g^{st}Ag.$$  One checks using Lemma \ref{factor} that this action again preserves the inner product on ${\mathbb R}^{2,1|2}$ and hence restricts to an isometric
action on $\hyperp$ itself.

\medskip

\noindent 
The body of this action of ${\rm OSp}(1\vert 2)$ on $\hyperp$ is the classical action of orientation-preserving isometries on $\hyperp'$, and this extension is our analogous action of this Lie super group on the super hyperbolic plane.   Moreover, the induced action of $\widehat{\rm SL}_2<{\rm OSp}(1|2)$ 
preserves the bosonic hyperboloid $\widehat{\hyperp}$, and we therefore find the equivariant tower of isometric actions
$$\begin{aligned}
&{\rm OSp}(1|2)&\acts&\hskip 5ex\hyperp\\
&\hskip 3ex\vee&&\hskip 5ex\cup\\
&\hskip 2.5ex\widehat{\rm SL}_2&\acts&\hskip 5ex\widehat{\hyperp}\\
&\hskip 3ex\vee&&\hskip 5ex\cup\\
&{\rm SL}(2,{\mathbb R})&\acts&\hskip 5ex\hyperp'\, .\\
\end{aligned}$$

\medskip

\section{${\rm OSp}(1|2)$-invariant Area Form on $\hyperp$}\label{sec:omega}

The classical hyperboloid $\hyperp'$ supports the ${\rm SL}(2,{\mathbb R})$-invariant hyperbolic area form given in coordinates by
$\Omega'(x_1,x_2,y)={{dx_2\wedge dx_1}\over {2y}}$, as we shall recall.  We promote this to an $\widehat{\rm SL}_2$-invariant form on $\widehat{\hyperp}$ and then suitably pull-back to derive an ${\rm OSp}(1|2)$-invariant two-form on $\hyperp$ itself.

\begin{lemma} $\widehat{\hyperp}$ supports the $\widehat{\rm SL}_2$-invariant two-form
$\hat\Omega={{{\rm d}x_2\wedge {\rm d}x_1}\over{2y}}$.
\end{lemma}

\begin{proof}We have $\begin{psmallmatrix}a&b&0\\c&d&0\\0&0&1\\\end{psmallmatrix}.(x_1,x_2,y|0,0)=(x_1',x_2',y'|0,0),$ where
$$\begin{aligned}
x_1'&=a^2x_1+2acy+c^2x_2,\\
x_2'&=b^2x_1+2bdy+d^2x_2,\\
y'&=abx_1+(ad+bc)y+cdx_2.\\
\end{aligned}$$
One computes directly that
$$\begin{aligned}
\hat\Omega'&={{ {\rm d}x_2'\wedge  {\rm d}x_1'}\over{2y'}}\\
&={{(ad-bc)[(ad+bc)(dx_2\wedge dx_1)+2 {\rm d}y\wedge (ab~ {\rm d}x_1-cd~ {\rm d}x_2)]}\over {2[abx_1+(ad+bc)y+cdx_2]}}.
\end{aligned}$$
Meanwhile, $x_1x_2-y^2=1$ implies that $ {\rm d}y={{x_1 {\rm d}x_2+x_2 {\rm d}x_1}\over {2y}}$.
Substituting this into the numerator finally yields
$\hat\Omega=\hat\Omega'$.
\end{proof}

\begin{lemma} For any ${\bf x}=(x_1,x_2,y|\phi,\psi)\in\hyperp$, consider the element
$u^{\bf x}=u(y\psi-x_2\phi,y\phi-x_1\psi)\in {\rm OSp}(1|2)$.  This transformation acts on ${\bf x}$ to give
$u^{\bf x}.{\bf x}=(1-\phi\psi)(x_1,x_2,y|0,0)\in\widehat{\hyperp},$  and this determines the canonical mapping $$\begin{aligned}\hyperp&\to\widehat{\hyperp}\\
(x_1,x_2,y|\phi,\psi)&\mapsto(\hat x_1,\hat x_2,\hat y|0,0)=(1+\phi\psi)(x_1,x_2,y|0,0).
\end{aligned}$$
\end{lemma}
\vskip -.3in\hfill\hbox{{\rm \qedsymbol}}

\bigskip

\begin{theorem} The hyperbolic area form on $\hyperp'$ pulls-back under the natural mapping
$\hyperp\to\hyperp'$ to the ${\rm OSp}(1|2)$-invariant two-form
on $\hyperp$ given by
$$\Omega(x_1,x_2,y|\phi,\psi)={{{\rm d}[(1+\phi\psi)x_2)]\wedge{\rm d}[(1+\phi\psi)x_1]}\over {2[(1+\phi\psi)y)]}}.$$
\end{theorem}

\begin{proof} This follows from the previous lemma and the factorization in Lemma \ref{factor}.
\end{proof}

\begin{lemma}\label{prim}
$\Omega={\rm d}\omega$ on $\hyperp$ for the primitive
$$\omega(x_1,x_2,y|\phi,\psi)={{y(1+\phi\psi)}\over 2}~{\rm d}{\rm log}_e\biggl ( {{x_1}\over{x_2}}\biggr).$$
\end{lemma}

\begin{proof}
Let $\hat y=y(1+\phi\psi)$ and likewise for $\hat x_1,\hat x_2$.  Notice that
$$
\hat x_1\hat x_2-\hat y^2=(1+2\phi\psi)(x_1x_2-y^2)
=(1+2\phi\psi)(1-2\phi\psi)
=1,
$$
so that
$$2\hat y\,{\rm d}\hat y=\hat x_1 {\rm d}\hat x_2+x_2 {\rm d}\hat x_1.$$
Thus,
$$\begin{aligned}
{\rm d}\omega&={\rm d}\biggl [{{\hat y}\over 2}~ ({{{\rm d}\hat x_1}\over{\hat x_1}} -{{{\rm d}\hat x_2}\over{\hat x_2}} )\biggr]\\
&={{\hat x_1{\rm d}\hat x_2+\hat x_2{\rm d}\hat x_1}\over {4\hat y}}~\wedge~
{1\over{\hat x_1\hat x_2}}[\hat x_2{\rm d}\hat x_1-\hat x_1{\rm d}\hat x_2]\\
&={{{\rm d}\hat x_2\wedge {\rm d}\hat x_1}\over {2\hat y}}\\
&=\Omega,\\
\end{aligned}$$
and finally
$\omega={{\hat y}\over 2}~{\rm d} {\rm log}_e ({{\hat x_1}\over{\hat x_2}})
={{\hat y}\over 2}~{\rm d} {\rm log}_e ({{x_1}\over{x_2}} )$.\end{proof}

\section{Line Integrals}\label{sec:line}

Let ${\bf P}=(p_1,p_2,p|\alpha,\beta)$ and ${\bf Q}=(q_1,q_2,q|\gamma,\delta)$ be distinct points in $\hyperp$, and define 
$$d=\langle{\bf P},{\bf Q}\rangle>1~{\rm and}~\ell=\sqrt{{d+1}\over {d-1}}.$$ It follows that 
$d={{\ell^2+1}\over{\ell^2-1}}$ and the distance between ${\bf P}$, and ${\bf Q}$ is given by $D={\rm cosh}^{-1} d$, with ${\rm exp}(D)=d+\sqrt{d^2-1}={{\ell+1}\over{\ell-1}}$.

By the proof of Corollary \ref{geocor}, the geodesic segment $\overline{\bf PQ}$ from ${\bf P}$ to ${\bf Q}$  is parametrized for $t={\rm exp}(s)$ by
$$
{\bf x}(t)={1\over 2}(t{\bf e}+t^{-1}{\bf f}),{\rm for}~1\leq t\leq{\rm exp}(D)={{\ell+1}\over{\ell -1}},\\
$$
where
$$\begin{aligned}
{\bf e}&=M~[(1-\ell){\bf P}+(1+\ell){\bf Q}],\\
{\bf f}&=N~[(1+\ell){\bf P}+(1-\ell){\bf Q}],\\
\end{aligned}$$
with
$$M={{\ell-1}\over{2\ell}}~{\rm and}~ N={{\ell+1}\over{2\ell}},$$
so that $MN={1\over {2(d+1)}}$.
According to Lemma \ref{prim}, the area primitive on $\overline{\bf PQ}$ is given in coordinates as
$$\omega(t)=\omega\circ{\bf x}(t)=[1+\phi(t)\psi(t)]~\bar\omega(t),$$
where
$$\bar\omega(t)={1\over 4}[tA+t^{-1}B]~d{\rm log}_e~{{tA_1+t^{-1}B_1}\over{tA_2+t^{-1}B_2}},$$
with 
$$\begin{aligned}
A=M[(1-\ell)p+(1+\ell)q]~&{\rm and}~B=N[(1+\ell)p+(1-\ell)q],\\
A_i=M[(1-\ell)p_i+(1+\ell)q_i]~&{\rm and}~B_i=N[(1+\ell)p_i+(1-\ell)q_i],\\
\end{aligned}$$
for $i=1,2$, and with
$$\begin{aligned}
2\phi(t)&=tM[(1-\ell)\alpha+(1+\ell)\gamma]+t^{-1}N[(1+\ell)\alpha+(1-\ell)\gamma],\\
2\psi(t)&=tM[(1-\ell)\beta+(1+\ell)\delta]+t^{-1}N[(1+\ell)\beta+(1-\ell)\delta].\\
\end{aligned}$$

\bigskip

It is useful here and hereafter to introduce the notation
$$S_i={{(q_i\alpha-p_i\gamma)(q_i\beta-p_i\delta)}\over {(qp_i-pq_i)^2}},~{\rm for}~i=1,2.$$

\medskip

\begin{lemma} For $i=1,2$, we have the identities\label{samesign}
$$\begin{aligned}
A_iB_i&={{(qp_i-pq_i)^2}\over{d^2-1}}(1-2S_i),\\
AB_i-BA_i&=2~{{qp_i-pq_i}\over{\sqrt{d^2-1}}}.\\
\end{aligned}$$
In particular, it follows that $A_i$ and $B_i$ have the same sign provided $(qp_i-pq_i)_\#\neq 0$.
In this case, if $p_i,q_i>0$, then also $A_i,B_i>0$, and
$p_i\leq q_i$ if and only if ${\rm exp}(D){{A_i}\over{B_i}}\leq1$.
\end{lemma}

\begin{proof} For the first identity, compute
$$\begin{aligned}
A_iB_i&=MN[((1-\ell)p_i+(1+\ell)q_i][(1+\ell)p_i+(1-\ell)q_i]\\
&={1\over{2(d+1)}}\biggl\{(1-\ell^2)(p_i^2+q_i^2)+2(1+\ell^2)p_iq_i\biggr \}\\
&={1\over{2(d+1)}}\biggl \{ \biggl (1-{{d+1}\over{d-1}}\biggr)(p_i^2+q_i^2)+2\biggl (1+{{d+1}\over{d-1}}\biggr)p_iq_i\biggr \}\\
&={1\over{d^2-1}}\biggl\{ -p_i^2-q_i^2+p_iq_i[p_1q_2+p_2q_1-2pq+2(\alpha\delta+\gamma\beta)]\biggr\}\\
&={1\over{d^2-1}}\biggl \{-p_i^2-q_i^2+p_i^2q_1q_2+q_i^2p_1p_2
-2p_iq_i[pq-(\alpha\delta+\gamma\beta)]\biggr\}\\
&={1\over{d^2-1}}\biggl\{ p_i^2(q^2-2\gamma\delta)+q_i^2(p^2-2\alpha\beta)-2p_iq_i[pq-(\alpha\delta+\gamma\beta)]\biggr\}\\
&={1\over{d^2-1}}\biggl\{ (qp_i-pq_i)^2-2(q_i\alpha-p_i\gamma)(q_i\beta-p_i\delta)\biggr\},\\
\end{aligned}$$
as was asserted.
The second identity follows directly from the definitions and elementary algebra.

For the final assertions, let us first note that
$$
{\rm exp}(D){{A_i}\over{B_i}}
={{(1-\ell)p_i+(1+\ell)q_i}\over{(1+\ell)p_i+(1-\ell)q_i}}
$$ since ${\rm exp}(D){M\over N}={{\ell+1}\over{\ell-1}}\,{M\over N}=1$.
Now, $0<q_i\leq p_i$ implies that $p_i+q_i> \ell(q_i-p_i)$, so $B_i>0$, whence $A_i>0$.  Furthermore,
${\rm exp}(D){{A_i}\over{B_i}}\leq 1$ is equivalent to ${(1-\ell)p_i+(1+\ell)q_i}\leq{(1+\ell)p_i+(1-\ell)q_i}$, which is in turn equivalent to $q_i\leq p_i$.
In the contrary case, $0<p_i\leq q_i$ implies that $p_i+q_i>\ell(p_i-q_i)$, so $A_i>0$, whence $B_i>0$.  Again, 
${\rm exp}(-D){{B_i}\over{A_i}}\leq 1$ is equivalent to ${(1+\ell)p_i+(1-\ell)q_i}\leq{(1-\ell)p_i+(1+\ell)q_i}$, which is in turn equivalent to $p_i\leq q_i$.
\end{proof}

Our aim is to explicitly compute $\int_{\overline{\bf PQ}} \omega$ in this section in order to apply Stokes' Theorem in the next section, and we first assert that
\begin{eqnarray}\bar\omega(t)={{{\rm d}t}\over 2}~\sum_{i=1}^2 (-1)^i~{{(A+t^{-2}B)B_i}\over{t^2A_i+B_i}}.\label {*}\end{eqnarray}
To see this, compute directly that
$$d{\rm log}_e~{{tA_1+t^{-1}B_1}\over{tA_2+t^{-1}B_2}}={{2t^{-1}(B_2A_1-B_1A_2)}\over{(tA_1+t^{-1}B_1)
(tA_2+t^{-1}B_2)}}$$
and by partial fractions that
$${1\over {(tA_1+t^{-1}B_1)
(tA_2+t^{-1}B_2)}}={1\over{B_2A_1-B_1A_2}}\biggl \{ 
{{B_2}\over{t^2A_2+B_2}}- {{B_1}\over{t^2A_1+B_1}}\biggr \}.$$
Eqn.\,\ref{*} follows immediately.

Next we compute directly from the definitions that
$$\begin{aligned}4\phi(t)\psi(t)&=(\alpha\delta+\gamma\beta)[(t^2M^2+t^{-2}N^2)(1-\ell^2)+2MN(1+\ell^2)]\\
&+\alpha\beta\,[tM(1-\ell)+t^{-1}N(1+\ell)]^2\\
&+\gamma\delta\,[tM(1+\ell)+t^{-1}N(1-\ell)]^2\\
\end{aligned}$$
Collecting like powers of $t$, we find
\begin{eqnarray}\label{**}
\phi(t)\psi(t)=Xt^2+Yt^{-2}+Z,
\end{eqnarray}
where 
$$\begin{aligned}
4X&=M^2[(1-\ell)^2\alpha\beta+(1+\ell)^2\gamma\delta+(1-\ell^2)(\alpha\delta+\gamma\beta)]\\
4Y&=N^2[(1+\ell)^2\alpha\beta+(1-\ell)^2\gamma\delta+(1-\ell^2)(\alpha\delta+\gamma\beta)]\\
4Z&=2MN[(1-\ell^2)(\alpha\beta+\gamma\delta)+(1+\ell^2)(\alpha\delta+\gamma\beta)]\,.\\
\end{aligned}$$


\begin{lemma} \label{lem:lion} Suppose that $0<p_1<q_1$ and $0<q_2<p_2$.  Then we have
$$\begin{aligned}
&\int_{\overline{\bf PQ}}\omega=(A_1B_2-B_1A_2)\biggl[ 
{{AX}\over{(\ell-1)A_1A_2}}-
{{BY}\over{(\ell+1)B_1B_2}}\biggr]\\\\
+&\sum_{i=1,2} (-1)^i(1-S_i)\biggl[1+Z-X{B_i\over A_i}-Y{A_i\over B_i}\biggr]\,{\rm arctan}\biggl[{{qp_i-pq_i}\over{p_i+q_i}}(1-S_i)\biggr]\\
&+\pi\,sign(qp_2-pq_2)(1-S_2)\biggl[1+Z-X{B_2\over A_2}-Y{A_2\over B_2}\biggr].\\
\end{aligned}$$
\end{lemma}

\begin{proof}
Setting
$$\begin{aligned}
W(t)&={{B_2}\over{t^2A_2+B_2}}-{{B_1}\over{t^2A_1+B_1}}\\
\end{aligned}$$
and combining Eqns.~\ref{*} and \ref{**}, we find that
$$2\omega=\biggl\{t^2AX+t^0[BX+A(1+Z)]+t^{-2}[AY+B(1+Z)]+t^{-4}BY\,\biggr\}\,W(t){\rm d}t.$$

\medskip

Together perturbing ${\bf P}$ and ${\bf Q}$ by ${\rm OSp}(1|2)$ if necessary to arrange that $(pq_i-qp_i)_\#\neq 0$, we have $A_i,B_i>0$ by Lemma \ref{samesign}
and our hypotheses that $p_i,q_i>0$, for $i=1,2$.
We may then apply the following indefinite integrals
$$\begin{aligned}
\int{{t^{2}dt}\over{t^2u+v}}&=-{\sqrt{v\over{u^3}}~{\rm arctan}\biggl({\sqrt{u\over v}~t\biggr)}}+{t\over{u}},\\
\int{{dt}\over{t^2u+v}}&={1\over {\sqrt{uv}}}~{{\rm arctan}\biggl({\sqrt{u\over v}~t\biggr)}},\\
\int{{t^{-2}dt}\over{t^2u+v}}&=-{\sqrt{u\over{v^3}}~{\rm arctan}\biggl({\sqrt{u\over v}~t\biggr)}}-{1\over{vt}},\\
\int{{t^{-4}dt}\over{t^2u+v}}&={\sqrt{u{^3}\over{v^5}}~~{\rm arctan}\biggl({\sqrt{u\over v}~t\biggr)}}-{{v-3ut^2}\over{3v^2t^3}}\\
\end{aligned}$$
to conclude that
$$\begin{aligned}
2\int_{\overline{\bf PQ}}\omega&=AXt\biggl [{{B_2}\over{A_2}}-{{B_1}\over{A_1}}\biggr]+BYt^{-1}\biggl [{{A_1}\over{B_1}}-{{A_2}\over{B_2}}\biggr]\\
&+{{AB_2-BA_2}\over{\sqrt{A_2B_2}}}\biggl[ 1+Z-X~{{B_2}\over{A_2}}-Y~{{A_2}\over{B_2}}\biggr]~{\rm arctan}\biggl(\sqrt{{{A_2}\over{B_2}}}~t\biggr)\\
&-{{AB_1-BA_1}\over{\sqrt{A_1B_1}}}\biggl[ 1+Z-X~{{B_1}\over{A_1}}-Y~{{A_1}\over{B_1}}\biggr]~{\rm arctan}\biggl(\sqrt{{{A_1}\over{B_1}}}~t\biggr)\\
\end{aligned}$$
evaluated between the limits $t=1$ and $t={\rm exp}(D)={{\ell+1}\over{\ell -1}}$.

\medskip

Now, we compute
$$\begin{aligned}
{\rm arctan}\biggl(\sqrt{A_i\over B_i}~t \biggr)\Biggr |_{t=1}^{t={\rm exp}(D)}
&={\rm arctan}\biggl(\sqrt{A_i\over B_i}~{\rm exp}(D)\biggr)-{\rm arctan}\biggl(\sqrt{A_i\over B_i}\biggr)\\
&={\rm arctan}\biggl\{{{\sqrt{A_i\over B_i}({\rm exp}(D)-1)}\over{1+{A_i\over B_i}{\rm exp}(D)}}\biggr\}+k_i\,\pi,\\
\end{aligned}$$
using standard formulas for the difference of arctangents, our hypotheses on $p_i,q_i$,  and the last part of Lemma \ref{samesign}, where $k_i=\begin{cases}0,&{\rm if}~i=1;\\1,&{\rm if}~i=2.\\\end{cases}$

Meanwhile, we have
$$\begin{aligned}
{{\sqrt{A_iB_i}\,[{\rm exp}(D)-1]}\over{A_i\,{\rm exp}(D)+B_i}}~{{\rm exp(-{D\over 2})}\over{\rm exp(-{D\over 2})}}
&={{\sqrt{A_iB_i}~{\rm sinh}({D\over 2})}\over{A_i\,{\rm exp}({D\over 2})+B_i\,{\rm exp}(-{D\over 2})}},\\
\\
\end{aligned}$$
and moreover,
$$\begin{aligned}
A_i\,{\rm exp}({D\over 2})+B_i\,{\rm exp}(-{D\over 2})&=\sqrt{{\ell+1}\over{\ell-1}}\,\biggl({{\ell-1}\over{2\ell}}\biggr)
\biggl[(1-\ell)p_i+(1+\ell)q_i\biggr]\\
&+\sqrt{{\ell-1}\over{\ell+1}}\,\biggl({{\ell+1}\over{2\ell}}\biggr)
\biggl[(1+\ell)p_i+(1-\ell)q_i\biggr]\\
&={{\sqrt{\ell^2-1}\over\ell}}\, (p_i+q_i)\\
\end{aligned}$$
since ${\rm exp}(D)={{\ell+1}\over{\ell-1}}$ implies that ${\rm cosh}({D\over 2})={\ell\over{\sqrt{\ell^2-1}}}$
and ${\rm sinh}({D\over 2})={1\over{\sqrt{\ell^2-1}}}$.

\bigskip

Finally using the first identity in Lemma \ref{samesign}, we have
$$\begin{aligned}
-k_i\,\pi+{\rm arctan}\biggl(\sqrt{A_i\over B_i}~t \biggr)\Biggr |_{t=1}^{t={\rm exp}(D)}&=
{\rm arctan}\biggl({{2\,\sqrt{A_iB_i}\,{\rm sinh}({d\over 2})\,{\rm cosh}({d\over 2})}\over{p_i+q_i}}\biggr)\\
&={\rm arctan}\biggl({{{\rm sinh}(D)}\over{\sqrt{d^2-1}}}\,{{|qp_i-pq_i|}\over{p_i+q_i}}\,(1-S_i)\biggr)\\
&={\rm arctan}\biggl({{|qp_i-pq_i|}\over{p_i+q_i}}\,(1-S_i)\biggr)\\
\end{aligned}$$
since $d={\rm cosh}(D)$.

\bigskip

It follows from Lemma \ref{samesign} that
$$
{{AB_i-BA_i}\over{\sqrt{A_iB_i}}}=2\, {{qp_i-pq_i}\over{|qp_i-pq_i|}}\,(1-S_i)=2\, {sign(qp_i-pq_i)}\,(1-S_i),
$$
and so
$$\begin{aligned}
&{{AB_i-BA_i}\over{\sqrt{A_iB_i}}}\,\Biggl[{\rm arctan}\biggl(\sqrt{A_i\over B_i}~t \biggr)\Biggr |_{t=1}^{t={\rm exp}(D)}\Biggr]\\\\
&=2\,(1-S_i)\,{\rm arctan}\biggl[{{qp_i-pq_i}\over{p_i+q_i}}(1-S_i)\biggr]\\
&+2\,\pi\,k_i\,{sign(qp_i-pq_i)}\,(1-S_i)\\
\end{aligned}$$
since arctangent is an odd function.
\end{proof}


\section{Area of a Super Triangle}\label{sec:area}

Let $\triangle$ denote a super triangle with vertices
$$\begin{aligned}
{\bf P}&=(p_1,p_2,p|\alpha,\beta),\\
{\bf Q}&=(q_1,q_2,q|\gamma,\delta),\\
{\bf R}&=(r_1,r_2,r|\varepsilon,\phi).\\
\end{aligned}$$

\begin{lemma}\label{ace}
Without loss of generality by applying an element of ${\rm OSp}(1|2)$ and possibly relabeling the vertices, we may assume
that we have $\alpha=\gamma=\varepsilon$, $0<p_1<q_1<r_1$ and $0<r_2<q_2<p_2$.
\end{lemma}

\begin{proof} First recall that the classical isometry from $\hyperp'$ to the upper half plane 
is given by $\Pi:(x_1,x_2,y)\mapsto{{i-y}\over{x_2}}$, where $i=\sqrt{-1}$.  (This extends to an analogous mapping
from $\hyperp$ to the super upper half plane with complex coordinates $z,\eta$ by $$(x_1,x_2,y|\phi,\psi)\mapsto z={{i-y-i\phi\psi}\over{x_2}},~~
\eta={{\psi}\over {x_2}}(1+iy),$$ but we shall not need this here.  Each of these mappings is appropriately equivariant.)

It is thus geometrically clear that by translation in $\widehat {\rm SL}_2$, we may arrange that $\Pi({\bf P}),\Pi({\bf Q}),\Pi({\bf R})$ each have
positive imaginary part and negligible squared real part divided by imaginary part.  This implies that $(p_1p_2), (q_1q_2), (r_1r_2)$ are all nearly
equal to unity.  Now perturb in $\widehat{\rm SL}_2$ so that these imaginary parts are all distinct and relabel if necessary so that $0<r_2<q_2<p_2$, from which it
then follows that also $0<p_1<q_1<r_1$.

Next further perturb the resulting $\triangle$ with an element of $\widehat{\rm SL}_2$ so that
$$t=(qp_1-pq_1)+(rq_1-qr_1)+(pr_1-rp_1)$$
has non-zero body.
Given fermions $\xi,\eta$, let us suppose that $u(\xi,\eta).{\bf P}=(p_1',p_2',p'|\alpha',\beta'),$
and likewise respectively $\gamma'$ and $\varepsilon '$ for ${\bf Q}$ and ${\bf R}$.
One computes that
$$\begin{aligned}
\alpha'&=p_1\xi+p\eta+\alpha(1+\hbox{$\tiny{3\over 2}$}\xi\eta),\\
\gamma'&=q_1\xi+q\eta+\gamma(1+\hbox{$\tiny{3\over 2}$}\xi\eta),\\
\varepsilon'&=r_1\xi+r\eta+\varepsilon(1+\hbox{$\tiny{3\over 2}$}\xi\eta).\\
\end{aligned}$$

We therefore seek fermions $\xi,\eta$ so that in particular
$$\begin{aligned}
0=\alpha'-\gamma'=(p_1-q_1)\xi\biggl[1+{3\over 2} {{(\gamma-\alpha)\eta}\over{p_1-q_1}}\biggr]+(p-q)\eta+(\alpha-\gamma).\\
\end{aligned}
$$
Multiplying through by $[1-{3\over 2}{{(\gamma-\alpha)\eta}\over{p_1-q_1}}]$
yields 
$$
(p_1-q_1)\xi+(p-q)\eta\biggl [1-{3\over 2}{{(\gamma-\alpha)\eta}\over{p_1-q_1}}\biggr]
=(\gamma-\alpha)\biggl [1-{3\over 2}{{(\gamma-\alpha)\eta}\over{p_1-q_1}}\biggr],$$
or equivalently
$(p_1-q_1)\xi+(p-q)\eta=\gamma-\alpha$.  We likewise find in particular for $\alpha'=\varepsilon'$ that $(p_1-r_1)\xi+(p-r)\eta=\varepsilon-\alpha$.

Thus, we find the necessary and sufficient linear system
$$\begin{pmatrix}p_1-q_1&p-q\\p_1-r_1&p-r\\\end{pmatrix} \begin{pmatrix}
\xi\\\eta\\\end{pmatrix}=\begin{pmatrix}\gamma-\alpha\\\varepsilon-\alpha\\\end{pmatrix}$$
with solution
$$\begin{aligned}
\xi&={1\over t}~\biggl[(p-r)(\gamma-\alpha)+(q-p)(\varepsilon-\alpha)\biggr],\\
\eta&={1\over t}~\biggl[(r_1-p_1)(\gamma-\alpha)+(p_1-q_1)(\varepsilon-\alpha)\biggr],\\
\end{aligned}$$
where the determinant
$$\begin{aligned}
t&=(p-r)(p_1-q_1)-(p-q)(p_1-r_1)\\
&=(qp_1-pq_1)+(rq_1-qr_1)+(pr_1-rp_1)\\
\end{aligned}$$
has non-zero body by construction.
\end{proof}

For the rest of this paper, we shall assume that the triangle $\triangle$ with vertices ${\bf P,Q,R}$ has been normalized and possibly relabeled
to arrange that $\alpha=\gamma=\varepsilon$, $0<p_1<q_1<r_1$ and $0<r_2<q_2<p_2$.  In fact, no subsequent fermionic product of degree two will include
$\beta\delta$, $\delta\phi$, or $\phi\beta$.  Thus $\alpha=\gamma=\varepsilon$ divides all subsequent fermionic
expressions, whence the product of any two of them must vanish.  We shall refer to this fact as being {\sl due to fermionic degree}.

\bigskip

Attempting  to mitigate the inevitable proliferation of notation, set
$$\begin{aligned}
\bar V&=(A_1B_2-B_1A_2)\biggl[ 
{{AX}\over{(\ell-1)A_1A_2}}-
{{BY}\over{(\ell+1)B_1B_2}}\biggr],\\
H_i&=(Z-X{B_i\over A_i}-Y{A_i\over B_i}),~{\rm for}~i=1,2,\\
\end{aligned}$$
in the notation of the previous section, so that
$$\begin{aligned}
\int_{\overline{\bf PQ}}\omega
&=\bar V+\pi\,sign(qp_2-pq_2)\,(1+H_2-S_2)\\
&+\sum_{i=1,2}(-1)^i(1+H_i-S_i)\,{\rm arctan}\biggl[{{qp_i-pq_i}\over{p_i+q_i}}(1-S_i)\biggr]\\
\end{aligned}$$
since $H_iS_i=0$ due to fermionic degree.
This equation, together with its analogues for $\overline{\bf QR}$ and $\overline{\bf RP}$,
provides a dreadful but explicit closed-form expression $\int_\triangle\Omega=\int_{\overline{\bf PQ}}\omega+
\int_{\overline{\bf QR}}\omega+\int_{\overline{\bf RP}}\omega$ for the area of $\triangle$
by Stokes' theorem.  

Rather than explicate this expression and with our eye on the Angle Defect Formula, let us instead proceed to compute ${\rm cos}(\int_{\overline{\bf PQ}}\omega)$
and ${\rm sin}(\int_{\overline{\bf PQ}}\omega)$.  
These are tractable since ${\rm cos}(a)=1$
and ${\rm sin}(a)=a$ for any super number $a$ with $a^2=0$, and since a product of any two factors
$\bar V, H_1, H_2, S_1, S_2$ must vanish due to fermionic degree.

\begin{theorem}\label{bars}
We have the identities
$$\begin{aligned}
{\rm cos}\,\int_{\overline{\bf PQ}}\omega={{C+\hat c}\over F}\,,~&~{\rm and}~~~{\rm sin}\,\int_{\overline{\bf PQ}}\omega={{S+ \hat s}\over F},\\ \end{aligned}$$
 where
 $$\begin{aligned}
& C=-2[(p+q)^2+(d+1)(1-pq)]~~{\rm and}~~S=(p+q)(p_1q_2-p_2q_1),\\
 &F= 2(d+1)\sqrt{p_1p_2q_1q_2},~~\hat c=SV-\bar c~~{\rm and}~~\hat s=-CV-\bar s\,,\\
 \end{aligned}$$
 with
 $$\begin{aligned}
 \bar c&=2[(q\alpha-p\gamma)(q\beta-p\delta)+(\alpha+\gamma)(\beta+\delta)]\\
 &\hskip 2ex+(S_1+S_2)(qp_1-pq_1)(qp_2-pq_2),\\
 \bar s&=S_1(p_2+q_2)(qp_1-pq_1)-S_2(p_1+q_1)(qp_2-pq_2),\\
 \end{aligned}$$
 and
 $$\begin{aligned}
 V=&~\bar V+\pi\,sign(qp_2-pq_2)(H_2-S_2)\\
&+\sum_{i=1,2}(-1)^i\,(H_i-S_i)~{\rm arctan}\,\biggl[{{qp_i-pq_i}\over{p_i+q_i}}(1-S_i)\biggr].\\
\end{aligned} $$
\end{theorem}

\begin{proof}
We have derived
$\int_{\overline{\bf PQ}}\omega=(U_2-U_1)+V+\pi\,sign(qp_2-pq_2)$
in Lemma \ref{lem:lion},
where
$$\begin{aligned}
U_i&={\rm arctan}\,\biggl[{{qp_i-pq_i}\over{p_i+q_i}}(1-S_i)\biggr],~{\rm for}~i=1,2.\\
\end{aligned}$$
Standard trigonometric identities give
 $$\begin{aligned}
- {\rm cos}\,\int_{\overline{\bf PQ}}\omega&={\rm cos}\,U_2~ {\rm cos}\,U_1+{\rm sin}\, U_2~{\rm sin}\, U_1\\
 &~~-V\,[{\rm sin}\,U_2~ {\rm cos}\,U_1-{\rm cos}\, U_2~{\rm sin}\, U_1],\\
 - {\rm sin}\,\int_{\overline{\bf PQ}}\omega&=[{\rm sin}\,U_2~ {\rm cos}\,U_1-{\rm cos}\, U_2~{\rm sin}\, U_1\\
 &~~+V\,[{\rm cos}\,U_2~ {\rm cos}\,U_1+{\rm sin}\, U_2~{\rm sin}\, U_1],\\
 \end{aligned}$$
where the leading minus signs derive from $\pi\,sign(qp_2-pq_2)$ and where
$V^2=0$ follows from considerations of fermionic degree, the latter of which implies that ${\rm cos}\, V=1$ and ${\rm sin}\, V=V$.

The usual identities
 ${\rm cos~arctan~x}={1\over\sqrt{1+x^2}}$, ${\rm sin~arctan~x}={x\over\sqrt{1+x^2}}$ imply that
$$\begin{aligned}
{\rm cos}\, U_i&={{p_i+q_i}\over{\sqrt{2\,p_iq_i\,(d+1)}}},\\
{\rm sin}\, U_i&={{(qp_i-pq_i)(1-S_i)}\over{\sqrt{2\,p_iq_i\,(d+1)}}},\\
\end{aligned}$$
for $i=1,2$, using the equation
$$
(p_i+q_i)^2+(qp_i-pq_i)^2=2[(d+1)p_iq_i+(qp_i-qp_i)^2\, S_i],
$$
which relies upon $\langle{\bf P},{\bf Q}\rangle=d$ and $\langle {\bf P},{\bf P}\rangle=1=\langle {\bf Q},{\bf Q}\rangle$.  A further similar calculation gives
$$\begin{aligned}
&(p_1+q_1)(p_2+q_2)+(qp_1-pq_1)(qp_2-pq_2)\\
&=2[(p+q)^2+(d+1)(1-pq)-(q\alpha-p\gamma)(q\beta-p\delta)-(\alpha+\gamma)(\beta+\delta)],
\end{aligned}$$
and a more elementary direct calculation gives
$$\begin{aligned}
(p+q)(q_1p_2-q_2p_1)&=(p_1+q_1)(qp_2-pq_2)-(p_2+q_2)(qp_1-pq_1).\\
\end{aligned}$$
It remains only to collect the fermionic terms.
\end{proof}


Adopting for $\overline{\bf QR}$ and $\overline{\bf RP}$ the analogous respective alphabetic notation  introduced for $\overline{\bf PQ}$, Corollary~\ref{bars}
provides the expressions
$$\begin{aligned}
{\rm cos}\,\int_{\overline{\bf QR}}\omega={{D+\hat d}\over G}\,~&~{\rm and}~~~{\rm sin}\,\int_{\overline{\bf QR}}\omega={{T+\hat t}\over G},\\\\
{\rm cos}\,\int_{\overline{\bf RP}}\omega={{E+\hat e}\over H}\,~&~{\rm and}~~~{\rm sin}\,\int_{\overline{\bf RP}}\omega={{U+ \hat u}\over H},\\
\end{aligned}$$
where for example
$$\begin{aligned}
&D=-2[(q+r)^2+(e+1)(1-qr)]~~{\rm and}~~T=(q+r)(q_1r_2-q_2r_1),\\
&E=-2[(r+p)^2+(f+1)(1-rp)]~~{\rm and}~~U=(r+p)(r_1p_2-r_2p_1),\\
&G= 2(e+1)\sqrt{q_1q_2r_1r_2}~~{\rm and}~~H= 2(f+1)\sqrt{p_1p_2r_1r_2},
\end{aligned}$$
with $e=\langle{\bf Q},{\bf R}\rangle$, $f=\langle{\bf R},{\bf P}\rangle$.

\begin{corollary}\label{P}We have the identity
$$\begin{aligned}
{\rm cos}\, \int_\triangle \Omega&=P+R,\\
\end{aligned}$$
where 
$$\begin{aligned}(FGH)\,P&=\hat c(DE-TU)-\hat s(DU+ET)\\
&+\hat d(CE-SU)-\hat t(CU+ES)\\
&+\hat e(CD-ST)-\hat u(CT+DS),\\
(FGH)\,R&=CDE-CTU-DSU-EST,\\
\end{aligned}$$
with $FGH=8(d+1)(e+1)(f+1)\,p_1p_2q_1q_2r_1r_2$.
\end{corollary}

\begin{proof}
This follows directly from Corollary \ref{bars} and its analogues for $\overline{\bf QR}$ and $\overline{\bf RP}$, using standard
trigonometric formulas. 
\end{proof}

\section{Angle Defect}\label{sec:defect}

\begin{lemma} \label{left}
The cosine of the sum of the interior angles of $\triangle$ is given by
$$L=1-{{1+2def-d^2-e^2-f^2}\over{(d+1)(e+1)(f+1)}}\,,$$
where we have $d=\langle{\bf P},{\bf Q}\rangle$, $e=\langle{\bf Q},{\bf R}\rangle$, and $f=\langle{\bf R},{\bf P}\rangle$ as before. 
\end{lemma}

\begin{proof} Let $x,y,z$ be the respective cosines of the interior angles of $\triangle$ at ${\bf P,Q,R}$, so for example
$$x=-\langle{{{\bf Q}-{\bf P}d}\over{\sqrt{d^2-1}}},{{{\bf R}-{\bf P}f}\over{\sqrt{f^2-1}}}\rangle={{df-e}\over{\sqrt{d^2-1}\sqrt{f^2-1}}}$$
by Corollary \ref{cor:tgt}.  The usual formula for the sum of arccosines gives
$$\begin{aligned}
{\rm arccos}\,x&+{\rm arccos}\,y+{\rm arccos}\,z\\
&={\rm arccos}\biggl [xyz-z\sqrt{(1-x^2)(1-y^2)}\\
&\hskip 8ex-y\sqrt{(1-x^2)(1-z^2}-x\sqrt{(1-y^2)(1-z^2)}\biggr],\\
\end{aligned}$$
with for example
$$
x\sqrt{(1-y^2)(1-z^2)}=(df-e)\, {{(1+2def-d^2-e^2-f^2)}\over{(d^2-1)(e^2-1)(f^2-1)}}\\
$$
by direct calculation, noting that $1+2def\geq d^2+e^2+f^2$ by the Hyperbolic Law of Cosines for
$\triangle$.

Thus, we find
$$\begin{aligned}
(d^2-1)&(e^2-1)(f^2-1)L=(ef-d)(df-e)(de-f)\\
&+(1+2def-d^2-e^2-f^2)(d+e+f-ef-df-ed).\\
\end{aligned}$$
Dividing this equation by the common factor
$(d-1)(e-1)(f-1)$ of its two sides finally yields the asserted identity.
\end{proof}

In the notation of the previous section, define
\begin{eqnarray}\label{R}
R={{CDE-CTU-DSU-EST}\over{FGH}}.
\end{eqnarray}
Let us also define
$I=p_1q_2+p_2q_1$,  $J=q_1r_2+q_2r_1$, and $K=r_1p_2+r_2p_1$,
with
$i=\alpha\delta+\gamma\beta=d-\mbox{${\small {1\over 2}}$}I+pq$,
$j=\gamma\varepsilon+\delta\phi=e-\mbox{${\small {1\over 2}}$}J+qr$ and
$k=\varepsilon\beta+\alpha\phi=f-\mbox{${\small {1\over 2}}$}K+rp$,
and with
$\Theta={{\alpha\beta}\over{p^2+1}}+{{\gamma\delta}\over{q^2+1}}+{{\varepsilon\phi}\over{r^2+1}}$.

\begin{lemma}\label{Q}We have $L+R=Q$ with
$$\begin{aligned}
8&(d+1)(e+1)(f+1)Q\\
=&-4\biggl[{{\varepsilon\phi I^2}\over{r^2+1}}  +  {{\alpha\beta J^2}\over{p^2+1}}  +  {{\gamma\delta K^2}\over{q^2+1}}\biggr]\\
&+8[Ii+Jj+Kk]-2[iJK+jIK+kIJ]\\
&+2\Theta[I^2(r^2-1)+J^2(p^2-1)+K^2(q^2-1)]\\
&-4\Theta[IJ(1+rp)+JK(1+pq)+IK(1+qr)]\\
&-8\Theta[I(r^2+1)(1-pq)+J(p^2+1)(1-qr)+K(q^2+1)(1-rp)]\\
&-2\biggl[
{{KJi(pq+1)}\over{(p^2+1)(q^2+1)}}+{{KIj(qr+1)}\over{(q^2+1)(r^2+1)}}+{{IJk(rp+1)}\over{(p^2+1)(r^2+1)}}
\biggr]\\\\
&+4(p+q)(q+r)(r+p)\\
&\biggl[
{{Ji+iJ}\over{(q^2+1)(r+p)}}  +  {{Jk+Kj}\over{(r^2+1)(p+q)}}  + {{Ki+iK}\over{(p^2+1)(q+r)}}
\biggr]\\\\
&+4{{(p+q)(q+r)(r+p)}\over{(p^2+1)(q^2+1)(r^2+1)}}\\
&\biggl[
Ii{{(1-pq)(r^2+1)}\over{p+q}}  +Jj{{(1-qr)(p^2+1)}\over{q+r}}  +Kk{{(1-rp)(q^2+1)}\over{r+p}}
\biggr]\\\\
&+8(p+q)(q+r)(r+p)\\
\Biggl[
&{{I\biggl({{\alpha\beta}\over{p^2+1}} + {{\gamma\delta}\over{q^2+1}}\biggr)-i}\over{p+q}}
+{{J\biggl({{\gamma\delta}\over{q^2+1}} + {{\varepsilon\phi}\over{r^2+1}}\biggr)-j}\over{q+r}}
+{{K\biggl({{\alpha\beta}\over{p^2+1}} + {{\varepsilon\phi}\over{r^2+1}}\biggr)-k}\over{r+p}}
\Biggr].\\\\
\end{aligned}$$
\end{lemma}

One critical conclusion, which is tantamount to the classical Angle Defect Theorem, is that the body $Q_\#=0$ vanishes.

\begin{proof} First note for example that
$$IJ=q_1q_2K+p_1q_2^2r_1+p_2q_1^2r_2,$$
so 
$$\begin{aligned}(p_1q_2-p_2q_1)(q_1r_2-q_2r_1)&=q_1q_2 K-p_1q_2^2-p_2q_1^2r_2\\
&=2q_1q_2K-IJ,\end{aligned}$$
and likewise
$$\begin{aligned}(q_1r_2-q_2r_1)(r_1p_2-r_2p_1)&=2r_1r_2I-JK,\\
(p_1q_2-p_2q_1)(r_1p_2-r_2p_1)&=2p_1p_2J-IK.\\
\end{aligned}$$
It follows that
$$\begin{aligned}
IJK&=q_1q_2K^2+(p_1q_2^2r_1+p_2q_1^2r_2)K\\
&=q_1q_2K^2+p_1p_2(q_1^2r_2^2+q_2^2r_1^2)+r_1r_2(p_1^2q_2^2+p_2^2q_1^2)\\
&=q_1q_2K^2+p_1p_2[(q_1r_2+q_2r_1)^2-2q_1q_2r_1r_2]\\
&+q_1q_2[(p_1q_2+p_2q_1)^2-2p_1p_2r_1r_2]\\
&=q_1q_2K^2+p_1p_2J^2+r_1r_2I^2-4p_1p_2q_1q_2r_1r_2,\\
\end{aligned}$$
and moreover that
$$\begin{aligned}
(p_1q_2-q_1p_2)(q_1r_2-q_2r_1)K&=q_1q_2K-p_1p_2J-r_1r_2I+4p_1p_2q_1q_2r_2r_2,\\
(p_1q_2-q_1p_2)(r_1p_2-r_2p_1)J&=p_1p_2J-q_1q_2K-r_1r_2I+4p_1p_2q_1q_2r_2r_2,\\
(q_1r_2-q_2r_1)(r_1p_2-r_2p_1)I&=r_1r_2I-p_1p_2J-q_1q_2K+4p_1p_2q_1q_2r_2r_2.\\
\end{aligned}$$

Turning first to $L$ and setting $L'=8(d+1)(e+1)(f+1)L$, we find from Lemma \ref{left} that
$$L'=-def+d+e+f+de+ef+df+d^2+e^2+f^2.$$ A moderate computation
from the definitions using the expression above for $IJK$ and collecting like terms, we find that
$$\begin{aligned}
&L'=4p_1p_2q_1q_2r_1r_2\\
&+8[p^2q^2r^2+p^2q^2+q^2r^2+r^2p^2+pqr(p+q+r)-pq-qr-rp]\\
&+I^2[2-r_1r_2]+4I[1-pqr^2-r(p+q)-2pq]\\
&+J^2[2-p_1p_2]+4J[1-p^2qr-p(q+r)-2qr]\\
&+K^2[2-q_1q_2]+4K[1-pq^2r-q(r+p)-2rp]\\
&+2[
IJ(1+rp)+JK(1+pq)+IK(1+qr)]\\
&-2iJK+4i[J+K+Jrp+Kqr]+8i[1-pqr^2+I-2pq-r(p+q)]\\
&-2jIK+4j[I+K+Irp+Kpq]+8j[1-p^2qr+J-2qr-p(r+q)]\\
&-2kIJ+4k[I+J+Iqr+Jpq]+8k[1-pq^2r+K-2rp-q(r+p)].\\
\end{aligned}$$

Turning now to the calculation of R and using the previous formulas, we find that
$$\begin{aligned}
&-CTU=2(q+r)(r+p)(p^2+1)(q^2+1)[2r_1r_2I-JK]\\
&+2i(q+r)(r+p)(1-pq)[2r_1r_2I-JK]\\
&+(q+r)(r+p)(1-pq)[r_1r_2I^2-q_1q_2K^2-p_1p_2J^2+4p_1p_2q_1q_2r_1r_2],\\\\
&-DSU=2(p+q)(r+p)(q^2+1)(r^2+1)[2p_1p_2J-IK]\\
&+2j(p+q)(r+p)(1-qr)[2p_1p_2J-IK]\\
&+(p+q)(r+p)(1-qr)[p_1p_2J^2-r_1r_2I^2-q_1q_2K^2+4p_1p_2q_1q_2r_1r_2],\\\\
&-EST=2(p+q)(q+r)(r^2+1)(P^2+1)[2q_1q_2K-IJ]\\
&+2k(p+q)(q+r)(1-rp)[2q_1q_2K-IJ]\\
&+(p+q)(q+r)(1-rp)[q_1q_2K^2-p_1p_2J^2-r_1r_2I^2+4p_1p_2q_1q_2r_1r_2],\\\\
&-CDE=8(p^2+1)^2(q^2+1)^2(r^2+1)^2\\
&+(1-pq)(1-qr)(1-rp)[IJK+2IJk+2JKi+2IKj]\\
&+2(p^2+1)(q^2+1)(1-qr)(1-rp)[JK+2kJ+2jK]\\
&+2(q^2+1)(r^2+1)(1-pq)(1-rp)[IK+2iK+2kI]\\
&+2(p^2+1)(r^2+1)(1-pq)(1-qr)[IJ+2iJ+2jI]\\
&+4(p^2+1)^2(q^2+1)(r^2+1)(1-qr)[J+2j]\\
&+4(p^2+1)(q^2+1)^2(r^2+1)(1-rp)[K+2k]\\
&+4(p^2+1)(q^2+1)(r^2+1)^2(1-pq)[I+2i].\\
\end{aligned}$$

It is now purely a matter of using the earlier expression for $IJK$ and collecting like terms to derive the
asserted expression for $Q$.  The vanishing of the body of $Q$ is equivalent to the classical Angle Defect Theorem
and gives an internal consistency check.
\end{proof}

Combining the previous Theorem with Corollary \ref{P}, we finally have:

\begin{theorem} {\rm [The Angle Defect Theorem for Super Triangles]}\label{mainthm}

\noindent The angle defect ${\mathcal D}$ minus the area ${\mathcal A}$
of the super triangle $\triangle$ is given
by the fermionic correction $${\mathcal D}-{\mathcal A}={{P+Q}\over\sqrt{1-R^2}},$$ where $P$ and $Q$ are respectively given explicitly in Corollary \ref{P}
and Lemma \ref{Q}, and $R$ is expressed in Eqn.~\ref{R}.  This fermionic correction is not identically zero.
\end{theorem}

\begin{proof} We have derived in the referenced results that
$${\rm cos}\,{\mathcal A}=R+P~~{\rm and}~~{\rm cos}\,{\mathcal S}=Q-R,$$
where ${\mathcal S}$ is the sum of the interior angles of $\triangle$.
The standard formulas ${\rm arccos}\,x={\pi\over 2}-\sum_{n=0}^\infty 
{{(2n)!}\over{2^{2n}(n!)^2}}\,{{x^{2n+1}}\over{2n+1}}$ and ${{{\rm d}}\over{{\rm d}x}}\,{\rm arccos}\,x=-{1\over\sqrt{1-x^2}}$ therefore give
$${\mathcal A}={\rm arccos}\,R-{P\over\sqrt{1-R^2}},~~{\rm and}~~{\mathcal S}=\pi-{\rm arccos}\,R-{Q\over\sqrt{1-R^2}}$$
since $P^2=Q^2=0$.  The asserted formula therefore follows since by definition ${\mathcal D}-{\mathcal A}=(\pi-{\mathcal S})-{\mathcal A}$.

It is a triviality that our fermionic correction  is not identically zero on super triangles.
Indeed, suppose ${\mathcal D}-{\mathcal A}$ does vanish for some super triangle $\triangle$, so in particular ${\mathcal A}_\#={\mathcal A}$
by the classical Angle Defect Theorem.
Modify $\triangle$ by translating its two fermionic coordinates by respective odd constants $\phi$ and $\psi$
with $\phi\psi\neq 0$.  The resulting defect ${\mathcal D}'={\mathcal D}$ is unchanged, and yet
the resulting area ${\mathcal A}'={\mathcal A}(1+\phi\psi)={\mathcal A}_\#(1+\phi\psi)$ according to the formula for the area primitive.
Thus, our fermionic correction ${\mathcal D}'-{\mathcal A}'={\mathcal A}_\#\,\phi\psi$ cannot then also vanish for the modified triangle.

\end{proof}

\appendix

\section* {Appendix: Ideal Super Triangles}\label{appendixA}

In this appendix, we show that the fermionic part of the area of a super ideal triangle diverges.
To this end given ${\bf e}',{\bf f}',{\bf g}'\in{\mathcal L}^+$, define
$$
{\bf e}=\sqrt{{2\langle{\bf f}',{\bf g}'\rangle}\over{\langle{\bf e}',{\bf f}'\rangle\,\langle{\bf g}',{\bf e}'\rangle}}~~ {\bf e}',~~
{\bf f}=\sqrt{{2\langle{\bf g}',{\bf e}'\rangle}\over{\langle{\bf e}',{\bf f}'\rangle\,\langle{\bf f}',{\bf g}'\rangle}}~~ {\bf f}',~~
{\bf g}=\sqrt{{2\langle{\bf e}',{\bf f}'\rangle}\over{\langle{\bf f}',{\bf f}'\rangle\,\langle{\bf g}',{\bf e}'\rangle}}~~ {\bf g}'\,,
$$
so that $\langle{\bf e},{\bf f}\rangle=\langle{\bf f},{\bf g}\rangle=\langle{\bf g},{\bf e}\rangle=2$, and introduce coordinates
$$\begin{aligned}
{\bf e}&=(p_1,p_2,p|\alpha,\beta),\\
{\bf f}&=(q_1,q_2,q|\gamma,\delta),\\
{\bf g}&=(r_1,r_2,r|\varepsilon,\phi).\\
\end{aligned}$$
These define respective super geodesics
$$\begin{aligned}
E&=\{\mbox{${1\over 2}$}(t{\bf f}+t^{-1}{\bf g})\in\hyperp:0<t<\infty\}\,,\\
F&=\{\mbox{${1\over 2}$}(t{\bf g}+t^{-1}{\bf e})\in\hyperp:0<t<\infty\}\,,\\
G&=\{\mbox{${1\over 2}$}(t{\bf e}+t^{-1}{\bf f})\in\hyperp:0<t<\infty\}\,,\\
\end{aligned}$$
as in Theorem \ref{uvthm}, bounding a super ideal triangle.

The ${\rm OSp}(1|2)$-invariant two-form $\Omega$ on $\hyperp$ admits the same primitive $\omega$ as in Section \ref{sec:omega}, and 
we begin by computing $\int_G\omega$.  The calculations
in Section \ref{sec:line} apply directly to conclude that
$$2\omega=\biggl \{t^2pX+t^0[qX+p(1+Z)]+t^{-2}[pY+q(1+Z)]+t^{-4}qY\biggr\}\,U(t)~{\rm d}t,$$
where $4X=\alpha\beta$, $4Y=\gamma\delta$, $4Z=\alpha\delta+\gamma\beta$, and
$U(t)={{q_2}\over{t^2p_2+q_2}}-{{q_1}\over{t^2p_1+q_1}}$, and moreover that
$$\begin{aligned}
2\int_G\omega&=pXt\biggl [{{q_2}\over{p_2}}-{{q_1}\over{p_1}}\biggr]+qYt^{-1}\biggl [{{p_1}\over{q_1}}-{{p_2}\over{q_2}}\biggr]\\
&+{{pq_2-qp_2}\over{\sqrt{p_2q_2}}}\biggl[ 1+Z-X~{{q_2}\over{p_2}}-Y~{{p_2}\over{q_2}}\biggr]~{\rm arctan}\biggl(\sqrt{{{p_2}\over{q_2}}}~t\biggr)\\
&-{{pq_1-qp_1}\over{\sqrt{p_1q_1}}}\biggl[ 1+Z-X~{{q_1}\over{p_1}}-Y~{{p_1}\over{q_1}}\biggr]~{\rm arctan}\biggl(\sqrt{{{p_1}\over{q_1}}}~t\biggr),\\\end{aligned}$$
but in this instance evaluated between the limits $t=0$ and $t=\infty$.

The latter two summands are bounded, and the former two summands give
$$pXt\biggl [{{q_2}\over{p_2}}-{{q_1}\over{p_1}}\biggr]+qYt^{-1}\biggl [{{p_1}\over{q_1}}-{{p_2}\over{q_2}}\biggr]
=(p_1q_2-p_2q_1)\biggl[{{Xt}\over p}+{Y\over {qt}}\biggr].$$
This expression diverges for the specified limits provided only that we have $p_1q_2\neq p_2q_1$.

There can be no cancellation from
$\int_E\omega$ and $\int_F\omega$, since these involve different fermionic products, and
it follows that the fermionic part of the area of an ideal triangle indeed diverges.  One can readily trace
through the remainder of the calculations in the main text to confirm that the body of
this area equals $\pi$.


\bigskip

\begin{thebibliography}{ABCD}

\bibitem{bryce}
DeWitt, B. {\it Supermanifolds} 2nd Edition, Cambridge Monographs on Mathematical Physics,
Cambridge University Press, Cambridge, UK (1992).

\bibitem{kessler}
Kessler, E. {\it Supergeometry, Super Riemann Surfaces and the Superconformal
Action Functional}, Lecture Notes in Mathematics {\bf 2230},
Springer-Nature (2019).

\bibitem{HPZ}
Huang Y., Penner, R., Zeitlin, A.
Super McShane identity, to appear
{\it Journal of Differential Geometry}.

\bibitem{Kac} Kac, V. Lie superagebras, {\it Advances in Mathematics} \textbf{26} (1977),  8--96.

\bibitem{lambert} Lambert, J.H., Theorie der Parallellinien (Unpublished Manuscript) (1766)


\bibitem{Manin} Manin, Yu.~I. {\it Topics in Noncommutative Geometry},  Princeton University Press, Princeton N.J. (1991).

\bibitem{French} Papadopoulos, A. and Th\'eret, G. {\it La th\'eorie des parall\`eles de Johann Heinrich Lambert}. Critical edition with French translation and mathematical and historical commentaries, ed. Blanchard, coll. Sciences dans l'Histoire, Paris, 214 p., (2014).

\bibitem {PP1}Papadopoulos, A. and Penner, R.  A characterization of pseudo-Anosov foliations
{\it  Pacific Journal of Mathematics}   {\bf 130} (1987),  359-377.


\bibitem {PP2} ---, Enumerating pseudo-Anosov conjugacy classes, 
{\it Pacific Journal of Math}  {\bf 142}  (1990), 159-173. 	

\bibitem {PP3} ---, The Weil-Petersson symplectic structure at Thurston's
boundary,  {\it Transactions of the American Math Society}    {\bf 335 }(1993), 891-904.

\bibitem {PP4} ---,
La forme symplectique de Weil-Petersson et le bord de 
Thurston de l'espace de Teichm\"uller, {\it Comptes Rendus Acad.\ Sci.\ Paris} {\bf 312} 
S\'erie I (1991), 871-874.

\bibitem {PP5}
---,
Broken hyperbolic
structures,  {\it Annals of Global Analysis and Geometry} {\bf 27} (2005), 53-77.


\bibitem {PP6}
---, Hyperbolic metrics, measured foliations
and pants decompositions for non-orientable surfaces, 
Asian Journal of Mathematics {\bf 20} (2016), 157-182.



\bibitem{pennerbook}
Penner, R. {\it Decorated Teichm\"uller Theory}, QGM Masters Class Series, volume 1, European
Math Society  (2012).

\bibitem{norbert}
---, Super Hyperbolic Law of Cosines: same formula with
different content, Festschrift for Norbert A'Campo's 80th birthday.

\bibitem {PZ} Penner, R. and Zeitlin, A. Decorated super-Teichm\"uller space,
{\it Journal of Differential Geometry} {\bf 111} (2019), 527--566.

\bibitem{Rogers}
Rogers, A. {\it Supermanifolds: Theory and Applications},
World Scientific, Singapore (2007).







\bibitem{German} St\"ackel, P.  and Engel, F.  Die Theorie der Parallellinien von Euklid bis auf Gauss, eine Urkundensammlung zur Vorgeschichte der nicht-euklidischen Geometrie. B. G. Teubner, Leipzig, 1895.


\end{thebibliography}
\end{document}